\documentclass[amssymb]{amsart}
\usepackage{amssymb}
\usepackage{latexsym}
\usepackage{xypic}
\title[$(\Z/2)^n$-equivariant $K$-theory]{The $RO(G)$-graded
coefficients of $(\Z/2)^n$-equivariant $K$-theory}
\author[P. Hu and I. Kriz]{P. Hu and I. Kriz}
\input xypic
\newtheorem{theorem}{Theorem}
\def\Proof{\medskip\noindent{\bf Proof: }}

\def\Z{\mathbb{Z}}

\def\C{\mathbb{C}}

\def\C{\mathbb{C}}

\def\Pi{\mathbb{P}^{\infty}}

\def\qed{\hfill$\square$\medskip}

\def\Zpk{\mathbb{Z}/p^{k}}
\def\Zpk1{\mathbb{Z}/p^{k-1}}

\newcommand{\rref}[1]{(\ref{#1})}

\newcommand{\beg}[2]{\begin{equation}\label{#1}#2\end{equation}}
\def\r{\rightarrow}

\def\sl2{\widetilde{SL_{2}(\Z)}}

\begin{document}
\maketitle


\section{Introduction}

In physics, in type IIA and IIB string theory, $D$-brane charges
are calculated by $K$-theory (see \cite{mm} or 
\cite{w} for a survey). Kriz, Pando and Quiroz \cite{kpq},
following up on previous investigation of
other authors (e.g. \cite{geberd,geberd1})
investigated the basic case of $D$-brane charges of orbifolds
obtained by linear operation of finite groups on flat spacetime.
In the case of the finite group $(\Z/2)^n$, the groups
relevant can be expressed as the $RO(G)$-graded coefficients
of $(\Z/2)^n$-equivariant $K$-theory. In case of what is known
in physics as discrete torsion (cf. \cite{geberd1}), we get 
twisted $(\Z/2)^n$-equivariant $K$-groups with compact supports
of representations. A physical phenomenon called $T$-duality
predicts further certain relations between these groups.

\vspace{3mm} 
This is an amusing problem, and we found a direct
elementary solution, which is the subject
of the present note. We show that $K_{(\Z/2)^n}(S^V)$
for any finite representation $V$ is always concentrated in 
dimensions of one parity (even or odd) essentially $RO((\Z/2)^m)$
for some $m\leq n$ (see Theorem \ref{t1}). The exponents are
calculated by an explicit combinatorial algorithm. The twisted groups
are the same with a certain shift, which corresponds to
the physical $T$-duality prediction.

\vspace{3mm}
Since originally posting this note, Max Karoubi pointed out to us that
Theorem 1 can in fact be deduced 
as an easy consequence of his much more general result,
namely Theorem 1.8 of \cite{kar1}. We thank him for bringing this
to our attention. We also thank John Greenlees for previous discussions.

\vspace{5mm}

\section{The Main Theorem}

For a group $G$, we consider the reduced $G$-equivariant $K$-theory
$\tilde{K}^{\ast}_G$. 
The goal of this note is to prove the following theorem: 

\begin{theorem}
\label{t1}
For any finite dimensional representation $V$ of $G=(\Z/2)^n$, there 
exists an $\epsilon \in \{ 0, 1 \}$ and $m  \in \{ 0, 1, 
\ldots, n \}$, such that 
\begin{equation}
\begin{split}
\tilde{K}^k_G(S^V) & = \Z^{2^m} \mathrm{where\ } k \equiv \epsilon \ 
\ \mathrm{mod\ } 2 \\
& = 0\ \  \mathrm{else} . 
\label{main}
\end{split}
\end{equation}
Here, the number $m$ can be calculated by induction. 
\end{theorem}

\Proof
The basic observation is that by equivariant Bott periodicity
\cite{a}, the $K$-theory 
group~(\ref{main}) does not change when we add a complex representation 
$W$ to $V$. Also, the result is obviously true when $V$ is a sum of $\ell$ 
irreducible real representations $\alpha_1,...,\alpha_{\ell}$, 
which are $\Z/2$-linearly independent in the 
character group. In that case, we have $m = n-\ell$ and 
$\epsilon \equiv 0$ by smashing the usual cofiber sequence 
\begin{equation}
G/Ker(\alpha_i)_{+} \rightarrow S^0 \rightarrow S^{\alpha_i}
\label{cofiber}
\end{equation}
for $i=1,...,\ell$. 

Next, we recall that we also have Bott periodicity for equivariant 
$Spin^c$-representations (\cite{a}). If $\alpha, \beta, \gamma$ are real 
irreducible, 
\beg{e*}{ U = 1 + \alpha + \beta + \gamma + \alpha \beta + \alpha \gamma
+ \beta \gamma + \alpha \beta \gamma }
is spin. To this end, note that if we set 
\[ H^{\ast}((\Z/2)^3, \Z/2) = \Z[x, y, z] \]
then the total Stiefel-Whitney class of $U$ is 
\[ (1+x)(1+y)(1+z)(1+x+y)(1+x+z)(1+y+z)(1+x+y+z)  \]
which has even coefficients in degree $<3$. Consequently,
$w_1(U) = w_2(U) = 0$. 
So we have proved 
\begin{equation}
\tilde{K}^k_G(S^U) \cong \tilde{K}^k_G(S^0) . 
\label{realcase}
\end{equation}

It is now time for us to set up our induction process. We will choose 
irreducible real representations $\alpha_1, \ldots, \alpha_p$, 
independent in the class group, such that $V$ is a sum of tensor 
products of subsets of $\alpha_1, \ldots, \alpha_p$. In fact, we shall
represent $V$ as a ``hypergraph'', i.~e. the system $\Gamma$ of subsets of 
$\{ \alpha_1, \ldots, \alpha_p \}$ whose sum is $V$ (assuming it contains 
no complex sub-representations). By~(\ref{realcase}), we may assume 
$\Gamma$ is in fact just a ``graph'', i.~e. it does not contain any sets 
of cardinality $>2$ (because otherwise we
can add representations of the  form \rref{e*}
and subtract complex representations to reduce the number of sets
of highest cardinality). Note that we do not, however, assume $\Gamma$ contains 
all its ``vertices'', i.~e. $1$-element sets. 

\vspace{3mm}
Let us now look at a vertex $v$. We will do induction with respect to the number of 
edges adjacent to $v$, or alternately the number of vertices. 
If $v$ is attached 
to at least two edges $a, b$ in $\Gamma$ which in turn
attached to vertices $u,w$ (whether $v,w,u$ are included 
in $\Gamma$ or not), 
apply \rref{e*} with $\alpha=v$, $\beta=w$, $\gamma=u$. 
Then the graph $\Gamma$ will turn into a hypergraph again, which will 
contain the set $\{u,v,w\}$ not
contain either $a$ or $b$, and may contain some of the sets $\{u\}$,
$\{v\}$, $\{w\}$, $\{u,w\}$ (in fact precisely those which were not
contained in $\Gamma$). Now apply a base change which uses all the
original basis elements (=vertices) with the exception of $w$, which
will be replaced by $w^{\prime}=w+u$. Then representing the representation
again as a hypergraph with respect to the new basis, the set $\{u,v,w\}$
will turn into $\{v,w^{\prime}\}$. We see that there are now now more
sets of cardinality $>2$ attached to $v$, and the number of edges attached
to $v$ has decreased by $1$. Any sets of cardinality $>2$ in the
new graph not containing $v$ can be eliminated again as before by
(possibly repeated) application of \rref{e*}. 

\vspace{3mm}
By this induction, it suffices to consider the case where $v$ is attached 
by a single edge $a$ in $\Gamma$ (whether $v$ is included in $\Gamma$ or not). 
But the case when $v$ is not included in $\Gamma$ in fact equivalent to the 
graph obtained
by erasing $v$. Suppose therefore that $v$ is included in $\Gamma$.
Then let $v$ be attached by an edge to a vertex $w$ (which may or may not
be included in $\Gamma$). Now the point is that the entire induction
described above may be also applied to $w$ instead of $v$ to eliminate
all edges originating in $w$, with the possible exception of one
single edge, ending in a vertex $u$ (which may or may not be included
in $\Gamma$). Let us recapitulate then: we know that $\Gamma$ contains
the sets $\{v\}$, $\{v,w\}$, $\{w,u\}$, possibly (but not necessarily)
the set $\{w\}$, and not any other sets containing $v$ or $w$.
Now apply \rref{e*} with $\alpha=v,\beta=w,\gamma=u$. Then the new
hypergraph will contain the sets $\{u,v,w\}$, $\{u,v\}$ , possibly 
(but not necessarily) the set $\{w\}$ (in fact, precisely when $\Gamma$
didn't contain it), but no other sets
containing $v$ or $w$. Then, change basis by including all the vertices of
$\Gamma$ with the exception of $v$, which will be replaced by $v^{\prime}=u+v$.
Then the new hypergraph will then contain the sets $\{v^{\prime}\}$,
$\{v^{\prime},w\}$ and possibly the set $\{w\}$, and no other sets containing
$v^{\prime}$ or $w$. We see then that the hypergraph is a disjoint union
of two hypergraphs, which reduces this case to lower rank cases by
the K\"{u}nneth theorem.

\vspace{3mm}
Perhaps surprisingly, there is still a case which remains to be treated, 
namely the case of $n=2$, $V=\alpha+\beta+\alpha\beta$ where $\alpha$,
$\beta$ are the generators of the character group. Let $A,B,C\cong \Z/2$
be the subgroups of $(\Z/2)^2$ on which the representations 
$\alpha$, $\beta$, $\alpha\beta$ vanish. Then consider the cofiber
sequence
\beg{e+}{(\Z/2)^2/C_+\wedge S^{\alpha+\beta}\r S^{\alpha+\beta}\r
S^{\alpha+\beta+\alpha\beta}.
}
The $K^{1}_{(\Z/2)^2}$-groups of the first two terms vanish,
theis $K^{0}_{(\Z/2)^2}$-groups are $\Z\oplus\Z$, $\Z$ respectively,
so we need to identify the map 
\beg{ei}{\Z\oplus\Z\leftarrow \Z
}
induced on $K^{0}_{(\Z/2)^2}$ by the first map \rref{e+}. Note that
the interesting information is just the image of this map, which can be
calculated in $C\cong\Z/2$-equivariant $K$-theory. When restricted to
$C$, $\alpha\cong\beta$ are isomorphic to the sign representation of
$C$. Now in $K^{0}_{C}(S^{\alpha})$, we have an element $c$ which,
under the inclusion
\beg{e++}{S^0\subset S^{\alpha}
}
restricts to
$$1-\alpha\in R(C)=\tilde{K}^{0}_{C}(S^0).$$
Then the restriction of the generator of $\tilde{K}^{0}_{(\Z/2)^2}(S^{\alpha
+\beta})$ to $\tilde{K}^{0}_{C}(S^{2\alpha})$ is $c^2$. Our question
is thus equivalent to finding the image of $c^2$ in $\tilde{K}^{0}_{C}(S^0)$
under Bott periodicity. Now the Bott element
$$u\in \tilde{K}^{0}_{C}(S^{2\alpha})$$
maps to $c$ (by the construction of the Bott element as $1-H$ where
$H$ is the tautological line bundle on $\C P^1$ where $\Z/2$ acts
as $-1$ on $\C\subset\C P^1$). So, we need to find the image
of $u^2\in\tilde{K}^{0}_{C}(S^{4\alpha})$ under the composition
$Bf^*$ where
$f:S^{2\alpha}\r S^{4\alpha}$
is the inclusion (all such inclusions are equivariantly homotopic) and
$B$ is Bott periodicity. But we already know
$$f^*(u^2)=u(1-\alpha)$$
so
$$Bf^*(u^2)=1-\alpha.$$
This generates a direct summand in the left hand side of \rref{ei}.
Hence, we conclude
$$\tilde{K}_{(\Z/2)^2}^{1}(S^V)=\Z,$$
$$\tilde{K}_{(\Z/2)^2}^{0}(S^V)=0.$$ 
\qed
\vspace{5mm}

\section{An example: $n=3$.}

To demonstrate the algorithm described in the proof of Theorem \ref{t1}
in action, let us consider $n=3$, with the character group generated
by representations $\alpha,\beta,\gamma$. Then the cases which
do not immediately reduce to $n=2$ are when the representation $V$
is related by an automorphism of $(\Z/2)^3$ by one of the following:
\beg{ee1}{\alpha+\beta+\gamma +\alpha\beta\gamma
}
\beg{ee2}{\alpha+\beta+\gamma +\alpha\beta+\beta\gamma
}
\beg{ee3}{\alpha+\beta+\gamma +\alpha\beta+\alpha\beta\gamma
}
\beg{ee4}{\alpha+\beta+\gamma +\alpha\beta+\alpha\gamma+\beta\gamma
}
\beg{ee5}{\alpha+\beta+\gamma +\alpha\beta+\beta\gamma+\alpha\beta\gamma
}
\beg{ee6}{\alpha+\beta+\gamma +\alpha\beta+\beta\gamma+\alpha\gamma
+\alpha\beta\gamma.
}
First note that in the case of \rref{ee6}, as noted above, $1+V$ is
spin, so in this case, 
$$m=3, \;\epsilon=1.$$
In the case \rref{ee1}, the algorithm will first convert the
hypergraph into a graph using \rref{e*}, but then immediately
add the same representation \rref{e*} again, to change to the basis
$\alpha$, $\beta^{\prime}=\beta\gamma$, $\gamma$. In this basis
$V=\alpha+\alpha\beta^{\prime}+\beta^{\prime}\gamma +\gamma$. The
algorithm then adds \rref{e*} again to give $\alpha\gamma+\beta^{\prime}
+\alpha\beta^{\prime}\gamma$, and changes base again to $\alpha$,
$\beta^{\prime}$, $\gamma^{\prime}=\gamma\alpha$, which gives
$\beta^{\prime}+\gamma^{\prime}+\beta^{\prime}\gamma^{\prime}$.
So, keeping track of dimensions, we get
$$m=2,\;\epsilon=0.$$
In the case \rref{ee2}, the algorithm adds \rref{e*} to get
$\alpha\gamma+\alpha\beta\gamma$, which under the base
change $\alpha$, $\beta$, $\gamma^{\prime}=\alpha\gamma$ becomes
$\gamma^{\prime}+\beta\gamma^{\prime}$. So the answer is
$$m=1,\;\epsilon=1.$$
For \rref{ee3}, the algorithm applies \rref{e*} to reduce
to the graph $\alpha\beta+\beta\gamma$ which are independent,
so 
$$m=1,\;\epsilon=1.$$
For \rref{ee4}, the algorithm applies \rref{e*} to get
$\alpha\beta\gamma$, and then base change to
$\alpha$, $\beta$, $\gamma^{\prime}=\beta\gamma$, 
to give $\alpha\gamma^{\prime}$, so the answer is
$$m=2,\;\epsilon=1.$$
FOr \rref{ee5}, again the algorithm reduces to a graph by
applying \rref{e*},
giving $\alpha\gamma$. So again,
$$m=2,\;\epsilon=1.$$

\vspace{5mm}

\section{The twisted case}

It is interesting to note that Theorem~\ref{main} also gives a 
calculation of all \emph{twisted} reduced $(\Z/2)^n$-equivariant 
$K$-groups of representations (see \cite{as} for
the definition of twisted $K$-theory). To see this, first note that since 
$S^V = V^c$ (the $1$-point compactification of $V$), reduced 
$K$-theory of $S^V$ is $K$-theory with compact support of $V$. 
But $K$-theory with compact support on a space $X$ is, by 
definition, twisted by the same group as $K$-theory for $X$
(since it is a direct limit of related $K$-theories of pairs 
$(X, U)$ where $U$ is compact). But $V$ is contractible, so the 
ordinary ``lower'' twistings (note that there are no 
higher twistings over a point) are classified by 
\begin{equation}
H^3_{Borel}(V, \Z) = H^3 ((\Z/2)^n, \Z) . 
\end{equation}
When we denote
\begin{equation}
H^{\ast}((\Z/2)^n, \Z/2) = \Z/2[x_1, \ldots, x_n] 
\label{generators}
\end{equation}
with the usual generators (coming from $H^1(\Z/2, \Z/2)$ of the 
factors), then in fact 
\begin{equation}
H^3((\Z/2)^n, Z) = \Lambda^2[x_1, \ldots, x_n ] . 
\label{exterior}
\end{equation}
Indeed, the Bockstein 
\[ \beta: H^2 ((\Z/2)^n, \Z/2) \rightarrow H^3 ((\Z/2)^n, Z) \]
is onto (the group is annihilated by $2$), and $\beta(x_i^2) =0$, 
which by~(\ref{generators}) gives the isomorphism~(\ref{exterior}). 

Now let $\alpha_1, \ldots, \alpha_n$ be the generators of the 
character group of $(\Z/2)^n$, so that 
\begin{equation}
w_1 (\alpha_i) = x_i .
\end{equation}
Then we note that
\begin{equation*}
\begin{split}
w(1+ \alpha_i + \alpha_j + \alpha_i \alpha_j) & = 
(1+x_i)(1+ x_j)(1+ x_i + x_j) \\
& = (1+x_i + x_j)^2 + x_i x_j \\
& = 1+ x_i^2 + x_j^2 + x_i x_j , 
\end{split}
\end{equation*}
so 
\begin{equation}
\beta w_2(1 + \alpha_i + \alpha_j + \alpha_i \alpha_j)  = 
\beta(x_i x_j),
\label{bockstein}
\end{equation}
which are the generators of~(\ref{exterior}). 

Now, however, the following fact holds in twisted $K$-theory. Suppose 
$X$ is a $G$-space, and $\eta$ is a $G$-equivariant even-dimensional 
orientable vector bundle over $X$ with total space $V$. Then $\eta$ 
induces a non-equivariant vector bundle on the 
Borel construction $EG \times_G X$, so we have equivariant 
Stiefel-Whitney classes in Borel cohomology 
\[ w_i(\eta) \in H^i_{Borel}(X, \Z/2) . \]
The class
\[ w_3 (\eta) = \beta w_2(\eta) \]
is the obstruction to the bundle $\eta$ being a $G$-equivariant $Spin^c$-bundle. 
Now there is a Bott periodicity isomorphism on twisted $K$-theory with compact 
support: 
\begin{equation}
K_{G, \tau}^{i, c}(V) \cong K^{i, c}_{G, \tau+ w_3(\eta)}(X) . 
\label{twistedbott}
\end{equation}
(The proof is just a straightforward modification of Atiyah-Bott's index 
argument to the twisted case; there is also a formal argument using 
parametrized equivariant spectra. A very similar argument
is made in \cite{dk}.) But by~(\ref{bockstein}), we see that 
the $w_3$-classes of the representations $1+ \alpha_i + \alpha_j + 
\alpha_i \alpha_j$ generate the twisting group. When a twisting has
\[ \tau = \sum \epsilon_{ij} \beta(x_i x_j) , \]
then by~(\ref{twistedbott}), therefore 
\[ \tilde{K}^i_{\tau}(S^V) = 
\tilde{K}^i (S^{V+ \sum \epsilon_{ij}(1+ \alpha_i + \alpha_j 
+ \alpha_i \alpha_j)}) \]
which reduces the twisted case to the untwisted. 


\end{document}